\documentclass[11pt]{article}
\onecolumn
\newtheorem{lem}{{\sc Lemma}}[section]              
\newtheorem{thm}{{\sc Theorem}}[section]            
\newtheorem{cor}{{\sc Corollary}}[section]          

\newcommand{\bc}[2]{\left(\!\!\begin{array}c#1\\#2\end{array}\!\!\right)}

\usepackage{latexsym} 

\newcommand{\e}{\equiv}

\title{ \Large\bf
Additive structure of $Z$(.) mod $m_k$ (squarefree)\\[1ex]
   and Goldbach's Conjecture}
\author{ {\sc Nico F. Benschop} -  {\it ~~Amspade Research}, ~~The Netherlands}
\date{}

\begin{document}
\maketitle

\begin{abstract}
The product $m_k$ of the first $k$ primes (2..$p_k$) has neighbours $m_k \pm 1$
with all prime divisors beyond $p_k$, implying there are infinitely many
primes [Euclid]. All primes between $p_k$ and $m_k$ are in the group $G_1(k)$
of units in semigroup $Z_{m_k}$(.) of mutiplication mod $m_k$. Due to its
squarefree modulus $Z_{m_k}$ is a disjoint union of $2^k$ groups, with as many
idempotents - one per divisor of $m_k$, which form a Boolean lattice $BL$.
The  $additive$ properties of $Z_{m_k}$ and its lattice are studied.
It is shown that each complementary pair in $BL$ adds to 1 mod $m_k$ ,
and each even idempotent $e$ in $BL$ has successor $e$+1 in $G_1$. It follows
that $G_1(k)+G_1(k) \equiv E(k)$, the set of even residues in $Z_{m_k}$, so
each even residue is the sum of two roots of unity, proving "Goldbach for Residues"
mod $m_k ~(GR)$. The prime units in $G_1(k)$ have principle (natural) values
in the corresponding set $G(k)$ of naturals $u<m_k$. A proof by contradiction
and finite reduction using epimorphism $G_1(k+1) \rightarrow G_1(k)$ mod $m_k$,
and verifying $GC$ for $4<2n<30 ~(k=3)$, yields a contradiction for $k$=3.
This establishes $GC$: Each $2n>4$ is the sum of two odd primes.
\\
The structure of $G_1(k)$ mod $m_k$ is illustrated by the next features. The
smallest unit in $G_1(k)$ is $p_{k+1}$, so its smallest composite is $(p_{k+1})^2$.
Hence its units between $p_{k+1}$ and $(p_{k+1})^2$ are all prime if considered as
naturals (their principle values in set $G(k)$), to be used as summands for successive
$2n < m_k$. For $k=3~(m_3$=30) it is shown by complete inspection that each
$2n$ with $4<2n<30$ is indeed the sum of two odd primes. For $k>3$ the addition
to obtain $2n < (p_{k+1})^2$ produces no carry, thus yielding a natural sum.
The known Bertrand Postulate: $p_{k+1} < 2p_k$, implies overlapping intervals
for successive $2n$, again yielding $GC$ by contradiction to $GR$.
\end{abstract}

\\ \hspace*{2cm} Boolean lattice, Goldbach conjecture.

Subject msc: ~11P32

\section{Introduction} 

Detailed analysis of the algebraic structure of modulo arithmetic is
persued, especially multiplication in relation to addition and exponentiation.
Addition and multiplication are associative operations, so semigroup structure
analysis provides a good perspective for basic problems in arithmetic [2,3,6,8]
such as {\it Goldbach's conjecture} of every even number $2n>4$ being the sum
of two odd primes.
The additive structure of multiplicative semigroups with squarefree moduli
is studied, in ring $Z(+, .)$ mod $m_k$. Choosing as modulus the product
$m_k$ of the first $k$ primes, all primes between $p_k$ and $m_k$ are in
the group of roots of 1 mod $m_k$, denoted as the group $G_1$ of units.
As shown (thm~\ref{thm:GR}), $G_1+G_1$ covers all even residues $2n$ in $Z$ mod $m_k$.

The direct product $Z_{rs} = Z_r \times Z_s$ of multiplications with
coprime component moduli $r$ and $s$, is represented by component-wise
multiplication [4]. Squarefree modulus $m_k$ implies $Z_{m_k}(.) =
Z_{p_1} \times~.~.~ \times Z_{p_k}$ is a direct product of multiplications
mod $p_i$. This direct product is analysed as an ordered {\it disjoint
union} of maximal $subgroups$ derived from the component semigroups
$Z_{p_i}$. The emphasis is on the $additive$ properties of idempotents,
and the "fine structure" of residue ring $Z(+,~.)$ mod $m_k$. Considering
the principle values of units in units group $G_1(k)$ transfer additive
results from residues to positive integers. So residues $u$ mod $m_k$ are
taken as naturals $u< m_k$, with upper bound $u+v < 2m_k$ so a
possible carry is at most 1, and for sum $u+v<m_k$ no carry is produced.
For instance in Lemma~\ref{lem:compl}: the sum of each pair of complementary
idempotents equals 1 mod $m_k$, yielding the natural sum $m_k+1$ for pairs
other than \{0,1\}.

{\bf Notation}:
The known number representation (base $m$) $n=c.m+r$ with carry $c$ and
rest $0 \leq r < m$ is used. Operation + is natural addition, which for
two summands $<m_k$ can produce a maximal carry of 1 (base $m_k$).
For residue arithmetic $c=0$. This in contrast to the usual
interpretation of a residue arithmetic closure as an image $Z/Z_m$ of the
integers $Z$, consisting of residue classes to express the irrelevance
of the carry in residue arithmetic (mod $m$). In short, these additive
and multiplicative interpretations, with residue values vs. residue classes,
correspond to $n=r+0.m$ resp. $n=r+Z.m$, resetting the carry to $0$ or to
set $Z$. The latter interpretation is a mix of numbers and sets, which is
cumbersome and not required in the present additive analysis. Furthermore,
$m$ will denote modulus $m_k$ if no confusion can arise, and ~$\equiv$~
denotes congruence mod $m_k$. Sections 4 and 5 interprete residues
$n$ mod $m_k$ as naturals $n < m_k$ by taking their principle value,
where addition is restricted to summands that produce no carry.

The idempotents ~$e^2 \equiv e$ of ~$Z_m$(.) play an essential role.
For prime modulus $p$ it is known that $Z_p$ has just two idempotents:
0 and 1 mod $p$. And all residues $1,~.~.~,p$-1, coprime to $p$, are in
some permutation generated as residues of powers $g^i$ of some
{\it primitive root} $g<p$ of unity [1]. They form an order $p-1$ cyclic
subgroup $G$ of $Z_p$, written $G = g^* \equiv \{g^i\} ~(i=1..p$-1), with
$g^{p-1} \equiv 1$. Hence $Z_p(.)$ is a cyclic group, adjoined to zero.

{\bf Summary}:
~The product $m_k$ of the first $k$ primes is used
for analysis of all primes and their additive properties. Each of the $2^k$
divisors $d$ of $m_k$ yields a maximal subgroup $G_d$ of $Z_{m_k}$ containing
all $n<m_k$ with the same set of prime divisors as $d$. The respective group
identities are the $2^k$ idempotents of $Z_{m_k}$, ordered as Boolean lattice
$BL$ ~[4][6]~ of which the additive properties are studied.

The additive properties of $Z_{m_k}$ are characterised by the successor $n+1$
of any $n$, especially of the idempotents. An essential additive property is
that each complementary pair of idempotents in $BL$ sums to 1 mod $m_k$ (lemma~\ref{lem:compl}),
and every even $e^2=e$ has successor $e+1$ in $G_1$, while $G_1+G_1$ covers
all $2n$ mod $m_k$. This residue version $GR$ of Goldbach's Conjecture ($GC$)
is extended, by considering the set of principle values (naturals) $G(k)$ of
the units in $G_1(k)$ for $k \ge 3$, to prove $GC$ for positive integers.
Results listed in the Conclusions may be new.

For completeness, these essential concepts [5][6] are reviewed in sections
1 and 2. Section 3 derives a 'Goldbach-for-Residues'($GR$) result. Sections
4 and 5 give the approach to Goldbach's conjecture, followed by conclusions.

\section{Lattice of groups}  

In modulus $m_k=\prod p_i ~(i=1 ~..~ k)$~ each prime factor has exponent one.
So $m_k$, having no square divisor, is called {\it square free}. The prime
divisors of $m_k$ are referred to as {\bf base primes}.

Residues $n$ with the same base-prime divisors as squarefree divisor $d~|~m_k$
form a {\it maximal subgroup} $G_d \subset Z_{m_k}$(.) with closure due to
all possible products having the same base primes. If $e$ is the identity
(idempotent) of $G_d$, then each $n$ in subgroup ~$G_d \e G_e$~ has a unique
$local$ inverse $n^{-1}$~ defined by ~$n.n^{-1} \e e$.

The $2^k$ divisors of $m_k$ correspond to as many subsets of the $k$
base primes. Each divisor $d$ of $m_k$ generates a finite cycle $d^*=\{d^i\}$
with an idempotent $\underline{d}$, the identity of subgroup $G_d$. Each
subgroup has just one idempotent as its identity. So ~$Z_{m_k}$ has $2^k$
disjoint subgroups $G_d$, one for each divisor $d$ of $m_k$, ordered in a
Boolean lattice as their identities are ordered, as follows.

\subsection{Ordering of commuting idempotents}  

$Z_{m_k}$ is a disjoint union of $2^k$ groups $G_d$, and the group identities,
the idempotents, form a {\it Boolean lattice}. In fact, {\it commuting
idempotents} ~$e^2=e,~f^2=f$ can be {\it ordered} $e \geq f$ whenever $ef=fe=f$,
in other words $e$ is identity for $f$. ~~This is readily verified to be an
ordering relation, being transitive, anti-symmetric and reflexive [4].

The lattice {\it meet} (greatest lower bound) operation  is modelled by
{\it multiplication}. The product of two commuting idempotents $e, f$ is
idempotent: $ef.ef = effe = efe = eef = ef$, while $e,f$ are left- and right-
identity for $ef$ since $e.ef = ef = fe= fe.e$, sothat $e \geq ef$, and
similarly $f \geq ef$. Also, $ef$ is the greatest idempotent ordered under $e$
and $f$, since ~$c \leq e$ and $c \leq f$ imply $c \leq ef$, which is easily
verified.

The {\it join} (least upper bound) of two idempotents is the idempotent
with the {\it intersection} of the corresponding baseprime sets. Idempotent
'1' at the top has the smallest base-prime set (empty), while '0' at the
bottom contains all base-primes since 0=$m$ mod $m$.

The sum of two idempotents is generally not an idempotent, nor is its
generated idempotent their lattice-join, except for complementary
idempotents, to be derived next.

\subsection{ Lattice of idempotents: add {\it vs} join }  

As shown earlier, the set of idempotents of $Z$ mod $m$ is closed under
multiplication, forming a lower semi-lattice [4,6]. {\it Multiplication} models
the {\bf meet} (glb: greatest lower bound) operation of two idempotents,
yielding an idempotent with the $union$ of the respective base-prime sets.

Notice that all primes $p: ~p_k<p<m_k$~ are 'units' in topgroup $G_1$. In the
base-prime set of any idempotent or subgroup they are considered equivalent
to 1 mod $m_k$. For instance, cycle 2* mod $m$ (in $G_2$) produces residues
$c.2^n$, where ~$c \in G_1$~ are relative prime to $m_k$, and $c$ has
prime divisors $p_r > p_k$. Residues in $G_1$ can occur as factor in
each $n \in Z_{m_k}$, according to their name of $units$ in $Z_{m_k}$.

The {\bf join} (least upper bound $lub$) of two idempotents follows by
$intersecting$ their baseprime sets, yielding an idempotent with their
common baseprimes. \\{\bf Def:} two idempotents $a,b$ are
{\it complementary} iff $ab \e 0$ and $lub(a,b) \e 1$.

The endomorphism '$.e$' ~for idempotents $e$ in commutative ~$Z_m(.)$~
models the lattice meet operation by multiplication, since for each
~$x,y \in Z_m: ~xy.e \e xy.e^2 \e xe.ye$ .\\
Although in general the sum of two idempotents is  not an idempotent,
the next exception is an essential additive property of $Z_m(.)$ :

\begin{lem} \label{lem:compl} 
~~For any squarefree $m>1$ with at least two prime divisors:\\
For each $complementary$ pair $\{a,b\} \neq \{0,1\}$ of idempotents
in $Z_m(.)$ holds~ $a+b=m+1$.
\end{lem}
\begin{proof}
The lattice of idempotents has order ~$2^k$, with~ $2^{k-1}$~ complementary
pairs. Consider a sublattice of order four: ${0,1}$ ~and any other
complementary pair ~${a,b}$. It must be shown that $a+b \e 1$ mod $m$.
Now idempotents ~$a,b$~ are complementary, ~so ~$ab \e 0$ mod $m$, implying :
\\$(a+b)^2 \e a^2+2ab+b^2 \e a+b$ (mod $m$), thus~ $a+b$~ is idempotent.
And ~$(a+b)a \e a^2+ba \e a$ mod $m_k$, ~so~ $a+b \geq a$, and similarly
~$a+b \geq b$. Hence $a+b \e 1$ mod $m$, because by $lub(a,b) \e 1$ the
only idempotent covering complementary $a$ and $b$ is 1. Clearly, for
$\{a,b\} \neq \{0,1\}$ holds $1<a+b<2m$ so $a+b=m+1$, with carry =1 (base $m$).
\end{proof}

In other words : complementary idempotents $a,~b$ have disjoint base-prime
sets $A$ and $B$, and union $A \cup B$ consists of all base-primes in
$m$. For square-free $m, ~a.b \e 0$ is the idempotent containing all base-primes.
And $join(a,b)$ has the trivial intersection $A \cap B =1$ as base-prime
set, relative prime to $m$, with corresponding idempotent '1' of $G_1$.

\begin{lem} \label{lem:auto} 
For squarefree modulus $m=2.odd: ~h=m/2$ is the lowest odd idempotent in
$Z_m(.)$ and~ $a \rightarrow a+h$ is the only additive automorphism of $Z_m(.)$
\end{lem}
\begin{proof}
Notice that $2h \e 0$, so for each even or odd pair $a,b$ in $Z_m$ holds $(a+b)h \e 0$.
Hence :\\ $(a+h)(b+h)\e ab +(a+b)h + h^2 \e ab+h$ , and only if $h^2 \e h$
this yields $a \rightarrow a+h$ as additive automorphism of $Z_m(.)$.
Furthermore, $h=m/2$ is the lowest odd idempotent, namely the image under
$+h$ of the lowest even idempotent 0 in $Z_m$ (for squarefree $m:$ no
divisors of 0 exist). It is readily verified that this morphism is
1-1 onto, mapping $Z_m(even)$ and $Z_m(odd)$ into each other.
\end{proof}

Now consider product $m= m_k = \prod_{i=1}^k p_i$ of the first $k$ primes.
Unit 1 is ordered at the top of the lattice of idempotents, being the
identity for all idempotents in $Z_m = \times_i ~~Z_{p_i}$. Top group
$G_1$ of all residues relative prime to $m$ misses all base primes.
Thus ~$G_1 = \times_i ~~C(p_i$-1) ~$[i=2..k]$ is a direct product
of $k-1$ cycles of periods $p_i -1$.

\begin{cor} \label{cor:G1=G2}
~In $Z$(.) mod $m$ with square-free $m= 2.odd$, and let $h=m/2$ then:\\ \hspace*{.5cm}
$Odd$ and $even$ top-groups are isomorphic ~$G_1 \cong G_2$ under additive
automorphism $+h$.\\
{\rm Note:} isomorphic max cycles~ $(2+h)^* \cong 2^*$ in $G_1$ and $G_2$
(e.g. $5<primes<25$ are $15 \pm 2^i$)
\end{cor}

\section{ Primes, composites and neighbours}  

{\bf Equivalent sum and difference :}
~$(-1)^2$=1 implies $-1 \in G_1$, ~so~ $G_1 \e -G_1$ ~hence :

{\bf (1)}~~~~~~~~~~~~$G_1 + G_1 \equiv G_1 - G_1$

So sums and differences of pairs in $G_1$ yield the same set of
residues mod $m$. ~~Notice that:\\~$(-n)^2=n^2$, ~so ~$n$ and $-n$
generate the same idempotent, thus are in the same subgroup:

{\bf (2)}~~~~For every group $G_d \subset Z_m :$
~if $n \in G_d$~ then so is ~$-n$, while ~$G_d+G_d \equiv G_d-G_d$.

{\bf Neighbours $n$+1 and $n$-1 in the lattice of $Z_m$ :}

For integers and residues: $~n$ and $n$+1 are coprime for each $n$ so their
prime divisors form disjoint sets. The same holds for $n$ and $n-1$. Then
one would expect $n$ and $n$+1 to be in complementary subgroups of $Z_m$.
More precisely, the subgroup ordering of their idempotents implies:

\begin{lem}~ 
For each $n \in Z_m$~ and base-prime complementary ~$\overline{n}:$~~
 $G_{n \pm 1} \geq G_{\overline{n}}$
\end{lem}
\begin{proof} Due to the subgroup ordering, a {\it subset} of baseprimes
disjoint from (complementary to) those in $n$ defines a subgroup ordered
above or equal to $G_{\overline{n}}$.
\end{proof}

Hence $e$+1 for any {\bf even} idempotent $e$ must be in an odd subgroup
$G_d$ that is ordered $G_d \geq G_{\overline{e}}$, with $\overline{e}$ the
complement of $e$ in the lattice of $Z_m$. In fact, as shown next: ~$e$+1
is in topgroup $G_1$.

\subsection{Each idempotent's successor is in $G_1$ or $G_2$}  

The sum of two complementary idempotents yields an idempotent namely 1
(lemma~\ref{lem:compl}), which is their join or least upper bound. This is an
exception, and in general idempotents do not sum to an idempotent, let alone
their join. For instance, in $Z_{10}$ with idempotents 1, 5, 6, 0 : $5+1=6$
is idempotent, but join(5,1)= 1. And join(6,1)= 1 while $6+1=7$ is not
idempotent, although 7 does generate the proper idempotent 1, due to:

\begin{lem} \label{lem:e+1}
~In $Z$(.) mod $m$, with square-free $m=2.odd$:\\
(a) \hspace{1in} Each {\bf even} idempotent $e$ has $e$+1 in $G_1$, and\\
(b) \hspace{1in} each {\bf odd} idempotent $d$ has $d$+1 in $G_2$. \\
(c) \hspace{1cm}
  For period $n$ of $e+1$ in $G_1$ mod $m_k$ holds: ~$e.(2^n-1) \equiv 0$.
\end{lem}
\begin{proof} (a,c):
~Given $e^2=e$, notice that ~$(e+1)(e-1) \equiv e^2-1 \equiv e-1$, so $e$+1
is identity for $e-1$, hence $G_{e+1} \geq G_{e-1}$ for $every$ idempotent $e$.
~Now $(e+1)^2 \equiv e^2+2e+1 \equiv 3e+1$, and in general expanding $(e+1)^n$,
with $e^i \equiv e$ for all $i>0$ and factoring out $e$, yields:

~~~~~~~$(e+1)^n \equiv 1 + \sum_{i=1}^n \bc{n}{i} e^i \equiv 1 + (2^n-1)e$

We need to show $c=(2^n-1)e \equiv 0$ for every even idempotent $e$, where
$n$ is the period of $e+1$, with corresponding odd idempotent $d=(e+1)^n=c+1$,
which equals 1 iff $c \equiv 0$. ~In fact it would suffice if $2^n-1$ is in
a group complementary to $G_e$ in the lattice of $Z_m$. The baseprimes in
$2^n-1$, which are all necessarily odd, would then complement those in
even idempotent $e$.

This can be seen as follows: $d^2=d$ implies $(c+1)^2 \equiv c+1$, hence
$c^2+c \equiv 0$,\\ so:~~~~
$(2^n-1)^2~e +(2^n-1)e \equiv (2^n-1)(2^n-1+1)e \equiv (2^n-1)2^n e \equiv 0$.
\\ Apparently, the odd baseprimes in $2^n-1$ complement at least those
in $e$ because their union is complete (product 0). This implies
$(2^n-1)e = c \equiv 0$, independent of the extra factor $2^n$. So :

{\bf (3)} ~~~~~$(e+1)^n \equiv 1+(2^n-1)e \equiv 1$, ~where $n$ is the
period of $e$+1 in $G_1$.

Part (b) is dual to (a), proven similarly by using $G_1 \cong G_2$ (lemma~\ref{lem:G1=G2})
\end{proof}

\begin{thm}~( {\it Goldbach for Residues} ~~$GR$ ): \label{thm:GR} \\   
\hspace*{.5cm} For squarefree $m_k= \prod p_i ~(i=1...k)$ with $p_1$=2,
 and $E$ the set of even residues mod $m_k$:\\
\hspace*{1cm}
  In $Z$ mod $m_k:~~E \equiv ~\{2n\}~ \equiv G_1+G_1 \equiv G_1-G_1$, ~so :
\\ \hspace*{.5cm} Each even residue in ~$Z_{m_k}$ is a sum or difference
  of two units.
\end{thm}
\begin{proof}
~In short write $G$ for $G_1$. Let $e$ be any even idempotent, then multiply
$e \in G-1$ (lem~\ref{lem:e+1}) on both sides by $G$. On the lefthand side this
yields $G.e = G_e$ which is the max-subgroup on $e$, and on the righthand side
$G(G-1)= G^2-G = G-G$, sothat $G_e \subseteq G-G$. ~Using (1) yields:
~$G_e \subseteq G-G=G+G$ for all even $G_e$, so $G+G$ covers all even
residues.
\end{proof}

This also holds for any even squarefree modulus $m=2.odd$. Theorem 3.1 can
be generalized to hold for naturals which are the principle values of the
units in groups $G_1(k)$, as shown next.

\section{Prime units and carry extension} \label{def:principle}           

{\bf Define} ~$G_1(k)$~ as group of units mod $m_k$, and the correponding
set $G(k)$ of principle values (naturals) $\{1,u\}$~ where $p_k<u<m_k$
with $u$ coprime to base primes $p \leq p_k$. Use set $P(k)$ of all primes
in $G(k)$. The emphasis in the sequel is on the principle (natural) values
in $G(k)$ of units in group $G_1(k)$.

The primes $p>p_k$ are congruent mod $m_k$ to units in $G_1(k)$, and all
those $p<m_{k+1}$ in $G_1(k+1)$ are covered by $G(k) + a~m_k ~(carry~a~:
~0 \leq a<p_{k+1})$. \\An example for $k$=3 with all units in $G_1(4)$ follows
(table 1). It illustrates the relation between the prime structures of $G_1(k)$
and $G_1(k+1)$, which is a generalization of the known fact that all primes are
congruent to $G_1(2)=\{1,5\}$ mod $m_2$=6: remove the numbers that have a base
prime as divider (re Eratosthenes' ~prime sieve). Here : $G_1(4) \cong G_1(3)$
mod $m_3$=30.

\begin{verbatim}
 Units:  1     7<   11    13    17    19    23    29 :   mod m3 = 30

 +30:   31    37    41    43    47   7^2<   53    59     p_{k+1} = 7

 +60:   61    67    71    73  7.11<   79    83    89     7.n (8x '<')
                                                         not in G(4)
 +90:  7.13<  97   101   103   107   109   113   7.17<  (baseprime 7)

+120: 11^2#  127   131  7.19<  137   139  11.13# 149     Composites #
                                                         Smallest 11^2
+150:  151   157  7.23<  163   167  13^2#  173   179      in G(4)

+180:  181 11.17#  191   193   197   199  7.29< 11.19#   mod m4= 210
\end{verbatim}
{\bf Table 1} Unit extensions: $G(k+1) = \{~u+a.m_k~\}~:~unit~u \in G(k),~carry~0<a<p_{k+1}$\\

The set of all units in $G_1(k+1)$ is generated as illustrated for $G_1(4)$ in
table 1, including all primes in $P(k+1)$. Each natural unit $u \in G(k)$
generates at most $p_{k+1}-1$~ primes ~$p=u+a.m_k \in P(k+1)$~, with
~$p_k< p <m_{k+1}$ and carry $0 \le a < p_{k+1}$. For large enough $2n$
there are several prime pairsums in $GC$ format (see diagonals of equal
carry-sum in table 2). Clearly $p_{k+1}$ is the smallest unit in $G_1(k)$,
so ~$(p_{k+1})^2$~ is its smallest composite, hence:\\

{\bf (4)} ~~All (natural) units $u$ with $p_{k+1} \leq u<(p_{k+1})^2$ in $G(k)$ are prime.

Moreover, principle values of composite units in $G_1(k)$ are necessarily
generated under multiplication by the corresponding prime principle values $>p_k$
of units in $G_1(k)$. The reverse process of unit reduction by multiples of $m_k$
yields the next lemma:

\begin{lem} \label{lem:GRepi}  
~~$G_1(k+1) \rightarrow G_1(k)$ mod $m_k$, symbolizes that group $G_1(k)$ is
an epimorphic image of $G_1(k+1)$ with $v=t+c.m_k$, relating each principle
value $t \in G_1(k)$ to $p_{k+1}$ principle values $v \in G_1(k+1)$ with
a carry $c$.
\end{lem}
\begin{proof}
The mappings $v-c.m_k \longrightarrow t$ form a morphism because
$v.w \equiv (t+c.m_k)(u+d.m_k) \equiv (t.u)+e.m_k \equiv t.u $ mod $m_k$,
where $e = (td+cu)+cd.m_k$.
\end{proof}

Notice that each natural $n<m_k$~ is represented uniquely by $k$ digits
of a multi base code using the successive baseprimes: $p_1~.~.~.~p_k$.
The $k-1$ lower significant digits are extended with a most
significant digit or carry $a<p_k$, of weight $m_k$.

This in contrast to the usual single base code, e.g. decimal, using
powers of ten.  The successive bases 2, 6, 30, 210, ... have maximal
digit values $p_k$-1:~ 1, 2, 4, 6, ... respectively. For instance decimal
$331=210+11^2=1.210+4.30+0.6+0.2+1$ ~yields 5-digit code 1 4 0 0 1.

All primes $p>$3 are congruent to \{1, 5\} mod 6, while primes $p>p_3=5$
are congruent to the eight prime residues \{1,7,~...~,23,29\} mod 30 in
$G_1(3)$, obtained from ~$G(2)=\{1,5\}$~ by $p_3$-1=5-1=4 extensions with
increment $m_2=6$, namely \{7, 11\} ; \{13, 17\} ; \{19, 23\} ; \{25, 29\}.\\

Composite $25=5^2$ is not coprime to 30, hence is not in $G(3)$. The other
seven extensions are all primes $p_3<p<30=m_3$, forming with 1 the 8 units
in $G_1(3)=C_2 \times C_4$, ~in fact of form $15 \pm 2^i$ (cor~\ref{cor:G1=G2}).
The 7-1=6 extensions $G(3)+a.m_3$ generate all 2.4.6 = 48 units in $G_1(4):$
the 5 composites with prime divisors $p>7$ (exclude eight non-units of form
7.n in table 1), identity 1 and all ~$48-6=42$~ primes in open interval ~$(7,210)$.

\subsection{Pair sums of carry extended units}   

{\bf Define} ~set $S_0(k)=G(k)+G(k)$ ~of pair sums of (natural) units.
\\Denote even numbers interval by set ~$E(k)=\{2n~|~4<2n<m_k\}$, and the set
of natural carry-extended units: $T_a(k)=G(k)+a.m_k~~(a<p_{k+1})$ in $G(k+1)$.
The set of baseprimes is extended by $p_{k+1}$, so its multiples in $G(k+1)$
are not units (see table 1 for expanding G(3) to G(4) by baseprime 7).
All other extended units of $G(k)$ are  units of $G(k+1)$ since none is
divisible by a baseprime $p \leq p_{k+1}$.

Table 2 shows these sums for $k$=2 and 3 (by commutation half an
array suffices). Notice that ~$G(2)=\{1,5\}$ with pair sums $S_0(2)=
\{2,6,10\}$, while pair sums $2n$ in $S_0(3)= G(3)+G(3)$~ covers all $2n$ with
$2p_4 \leq 2n<m_3=30$, where ~$G(3)=\{1,7,~.~.~,29\}$ coprime to ~$2.3.5=30=m_3$.
For 6, 8, 10 use 3 and 5 to avoid non-prime 1.
In fact all ~$2n>16$~ have several ~$GC$ pair sums, e.g. each ~$2n$~ in
~$S_0(2)+6c=\{2,6,10\}+6c$~ for $1<c<p_3=5$~ has distinct unit pair sums,
all of which are prime pair sums.

\subsection{Pair sums of primes in ~$G(3)$}  

{\bf Define} ~Set ~$S_{a+b}(k) = S_0(k)+(a+b)m_k$ of pairsums in
$T_a(k)+T_b(k)$ of extended (natural) units in $G(k+1)$, except multiples
of $p_{k+1}$ (e.g. $5^2 \notin G(3)$), with carrysum ~$0 \leq a+b < p_{k+1}$.

{\bf Table 2} Extension sums:
carry sum diagonals ~$a+b=c<5$~ cover $2n$ in $E(3)$ by $S_0(2)+6c$
\begin{verbatim}
 Ta+Tb |    0        1        2        3        4    : carry b (wgt 6)
_______#__1___5 #  7__11   13__17   19__23   25__29  : translations Tb
     1 | .2.  6 |<- 6= 3+3                   xx
 0   5 |  6 .10.|<- S0(2)= {2,6,10}                25=5^2: no unit in G(3)
       #--------#   S0(3)= {2,6,..2n..,28)          6=3+3, 8=3+5, 10=5+5
     7 |  8  12  .14. 18
 1  11 | 12 >16<  18 .22.
 a  Ta  ........ ........+--------*
    13 | 14  18   20  24 |.26. 30 | <- S4(2)= S0(2)+4.6
 2  17 | 18  22   24  28 | 30 .34.|         ={26,30,34}
                +--------*--------+
    19 | 20  24 | 26  30 |
 3  23 | 24  28 | 30  34 |
       +--------*--------+
   5.5 x 26  30 |
 4  29 | 30  34 |
 ------*--------+
 \end{verbatim}

For instance extend ~$m_2$=6 ~to~ $m_3$=30 ~($p_3=5$) then
translations ~$S_{a+b}$~ of ~$S_0(2)$=\{2,6,10\} yield ~$5-1=4$
diagonals ~of~ 2 $\times$ 2 sums with carries ~$a+b<p_3$~ (table 2):

~~~~$S_1\{8,12,16\} ,~~S_2\{14,18,22\} ,~~S_3\{20,24,28\} ,~~S_4\{26,30,34\}$

Extending ~$G(2)=\{1,5\}$ ~yields~ $G(3)=\{~G(2)+6a~|~0<a<5~\}$,
containing prime set ~$P(3)=\{15 \pm 2^i,~29\}<m_3$ ~($i$=1,2,3)
~where $5^2$ is not coprime to 30, so not in $G(3)$.

\begin{lem}  \label{lem:base} 
~(basis $k$=3): Let $P'(3)=P(3)\cup \{3,5\}$, then ~$P'(3)+P'(3)$ ~covers ~$E(3)$.
\end{lem}
\begin{proof}
By complete inspection. Increments 4 in ~$S_0(2)=\{2,6,10\}$~ cause
successive $S_c(2)=S_0(2)+c.m_2$ with carry increment $m_2$=6 to
$interlace$ for $2n<m_3$. Exclude non-prime 1 by including primes 3 and 5.
Then $E(3)=\{2n \in [~6,30)~\}$~ is covered by pair sums of primes
$p<m_3$ in $G(3)$, extended with primes 3,~5 in $G(1) \cup G(2)$.
\end{proof}

So pair sum set $S_0(3)$, adapted for the interlacing edge-effect by
including $\{3,5\}$, covers adjacent $2n$ in $E(3)$. Hence interlacing
does not occur for $k>$3, and a unique carry sum $a+b=c$ suffices for
covering successive $2n$ by unit pair sums, in adjacent and disjoint
extension sum intervals $S_c(k)$, while:

~~~~Each $2n$ in $E(k+1)$ ~has a unique carry sum c with
  $0 \leq c < p_k$, such that $2n \in S_c(k)$.

This is to be used as basis for $k>3$, ~first for unit pair sum sets ~$S_c(k)$.

{\bf Define:} The set $S_0(k)=G(k)+G(k)$ of pair sums of principle values of
units is called {\it complete} if it covers $E(k)$, otherwise it is $incomplete$.

\begin{lem} \label{lem:Sc(k)}          
~For ~$k \geq 3$:

~~~~Extended pairsum sets $S_c(k)$ ~for~ $0 \leq c < p_k$
       ~partition $E(k+1)$ ~iff~ $S_0(k)$ covers $E(k)$.
\end{lem}
\begin{proof}
 ~Extension sets $S_c(k)=S_0(k)+c.m_k$ are disjoint for different carries
$c< p_{k+1}$, and $\{x,y\}$ in distinct extension sum sets remain so under
any shift~ $s=c.m_k:~~x \neq y \iff x+s \neq y+s$. For distinct carrysums
$c<c'$ with $c'-c=d:~ S_c(k) ~\cap~ S_c'(k) = S_c(k) ~\cap~ (S_c(k)+d.m_{k-1})=
\emptyset$. Their union covers $E(k) ~only$ if pair sums $S_0(k)=G(k)+G(k)$~
cover $E(k)$. Because some $2n$ missing from $S_0(k)$ implies its translations
$2n'=2n+c.m_k$~ are also missing from all $S_c(k)$ with $c>0$.
\end{proof}

\subsection{Excluding composites in $G(k)$, baseprimes and 1 as summands}

The set $G(k)$ of principle values (naturals) of units in group $G_1(k)$
coprime to baseprimes $2~..~p_k$, contains $p_{k+1}$ as smallest prime,
so the smallest composite in $G(k)$ is $(p_{k+1})^2$. Notice that $G(3)$
has no composites since $(p_4)^2=49>30=m_3$. Furthermore, the natural units
$u \in G(4)$ are in interval $(7<u<210)$ with smallest prime $p_5=11$,
hence minimal composite $11^2=121$, so all units of $G(4)$ in $[11,11^2)$
(coprime to 2.3.5.7=210) are prime. By inspection all $2n$ in interval
[22~..~222] are covered by prime pair sums, of which those $2n<210$
involve no carry.

The known Bertrand's Postulate is useful (Chebyshev 1850, simplified
by S.Pillai 1944) to prove a complete cover of even naturals:

{\bf BP} ~{\it (Bertrand's Postulate)}:
   ~~For each $n>1$ there is at least one prime between $n$ and $2n$.

Notice that  Pillai's proof [7] has an induction base of $2n \leq 60$
(see present lemma~\ref{lem:base}). In order to guarantee prime summands,
consider only pair sums of units $u<(p_{k+1})^2$, the smallest composite
in $G(k)$. In fact using $p_{k+1}<2p_k$ by Bertrands Postulate (BP),
the smaller interval $p_k < u< 2(p_{k+1})$ already suffices.
Successive $k$ yield $2n$ in overlapping intervals by $BP$, thus
covering all $2n$ beyond the induction base $k$=3. The next lemma is
readily verified, regarding the absence of a carry for $k > 3$.

\begin{lem} \label{lem:nocarry}  
For (natural) units in $G(k)$ and prime pairsums 2n in
$2~p_{k+1} \leq 2n < (p_{k+1})^2$ : no carry is produced for $k \geq 4$
since sum upperbound $(p_{k+1})^2 < m_k ~(~(p_{4+1})^2 =121 < m_4=210$).
\end{lem}

Notice that for initial $G(2)=\{1,5\}$ mod 6 (table 2) the baseprimes
2 and 3 are not used in pair sum residues $G(2)+G(2)= \{2,6,10\}$.
Considering $2n>4$ (re Goldbach's conjecture): non-prime 1 is avoided
by $6=3+3$ and $8=5+3$, the only $2n$ requiring 3. Moreover, 12=5+7
and 16=5+11 are the only extension pair sums $<30$ with one summand
of carry=0, thus requiring baseprime 5 of $G(3)$.

\section{Proving $GC$ by induction, or by reduction and contradiction}   

{\bf Approach} :~
Consider $G(k)$ (sect.~\ref{def:principle}) as set of 'natural units' $<m_k$
congruent to the units in group $G_1(k)$ of residues mod $m_k$, as defined in
the previous section. In other words, consider only the principle values in $G(k)$
\index{principle values} of the residue units in $G_1(k)$. Let $2n$ be small
enough, namely $2~p_{k+1} \leq 2n < (p_{k+1})^2$ with necessarily only prime
unit summands. There are two ways of proving $GC$: either a direct proof by
induction over $k \ge 3$ of primesums $2n<(p_{k+1})^2$ in $G(k)+G(k)$, or
an indirect proof by finite reduction of $GR(k)$ (thm~\ref{thm:GR}) and
contradiction to $S_0(3)$ (lemma~\ref{lem:base}). A direct proof, restricting
$GR(k)$ to the following primepair sums (item 3), would consist of the next
five steps:\\[1ex]
1. As summands use the principle value set $G(k)$ of the prime units in $G_1(k)$.\\
2. Complete inspection of $GC$ for $k=3$ (lemma~\ref{lem:base}).\\
3. Primesums $2n<(p_{k+1})^2$ in $G(k)+G(k)$ for $k>$3 yield no carry (lemma~\ref{lem:nocarry}).\\
4. Successive such restricted intervals of $2n$ intersect (Bertrand's Postulate).\\
5. The union over $k \ge 3$ of such primesums $2n$ in $G(k)+G(k)$ yield $GC$.

Using theorem~\ref{thm:GR} ($GR$) and lemma~\ref{lem:GRepi} an indirect proof
by \index{reduction} finite reduction: $S_0(k+1) \longrightarrow S_0(k)$ and
contradiction to lemma~\ref{lem:base} [~$S_0(3)$ covers $E(3)$~] runs as follows.

\begin{thm}~(Goldbach's Conjecture)~~
    Each $2n>4$ is a sum of two odd primes. \label{thm:GC}
\end{thm}
\begin{proof}
By inspection (lemma~\ref{lem:base}) $GC$ holds for $4<2n<m_3=30$ ($k$=3).
Now assume $GC$ to fail for an even primepair sum $2n<(p_{k+1})^2$ in
$S_c(k) \subset G(k)+G(k)$, to guarentee prime summands (by eqn(4)).
Then lemma~\ref{lem:GRepi}: $G(k+1) \rightarrow G(k)$ mod $m_k$
implies $2n-c.m_k$ to be missing as primepair sum from $S_0(k)$, making it
incomplete by lemma~\ref{lem:Sc(k)}. This in turn, by $G(k) \rightarrow G(k-1)$
mod $m_{k-1}$, reduces to incomplete $S_0(k-1)$, etcetera, down to incomplete
$S_0(3)$, contradicting lemma~\ref{lem:base}. Notice that $S_0(3)$ contains
only primepair sums, by eqn(4). Combining this with the overlap (by Bertrand's
Postulate) of prime summand intervals for successive $k$, hence also of pairsum
intervals $2p_{k+1} \leq 2n<(p_{k+1})^2$ for $k \geq 3$, it follows that Goldbach's
Conjecture holds. \end{proof}

Regarding the values of prime summands that suffice to cover all even
naturals, the following can be said. Notice that in table 2 (for $2n<30$)
only primes $p \geq p_k$ are required to represent $2n \geq 2p_k$ in most
cases. However, exeptions occur if prime gap $p_{k+1}-p_k>2$. Then $2n+2$
requires a prime $p_{~k-1}$ or smaller: the larger the gap the smaller
$p_{~k-i}$ is required. See for instance (table 2): $2n = 2p_k +2=
16,~28,~40$ ~for $p_k=7,~13,~19$ respectively, which require $p_{k-1}$
as Goldbach summand, due to a gap $p_{~k+1} - p_k=4$ (versus gap 2 in
cases $p_k=5,~11,~17$).

\section{Conclusions}

Balanced analysis of multiplication and addition in relation to each
other, with finite square-free moduli $2...p_k$ yields a fruitful
analysis of prime sums (Goldbach), similar to that with prime power
moduli mod $p^k$ for $p$-th power sums (Fermat [3], Waring [8]).
In both  approaches the careful extension of residues with a $carry$
is essential for transferring additive structural results to integers.
This 'residue-and-carry' method, as used for proving FLT [3] and Goldbach's
Conjecture, is based on unique number representation by residue and carry:
using the associative (semigroup) properties of the residue closure,
combined with an induction proof by carry extension. As such it could
well serve as a generic method to solve other hard problems in
elementary number theory [6].

In fact, the semigroup $Z_m$(.) of multiplication mod $m$ is formed by the
{\it endomorphisms} of the additive cyclic group $Z_m$(+) generated by 1.
 So ~$Z_m(.)=endo[~Z_m(+)~]$ where (.) distributes over (+), suggesting
a strong link between these two operations, evident from the derived
additive fine structure of $Z_{m_k}$ for squarefree modulus $m_k$.
A two-dimensional table of prime pair sums revealed additive properties
of ~$2n<m_3=30$ as basis for the analysis, hard to find otherwise.

The product $m_k$ of the first $k$ primes as modulus restricts all
primes between $p_k$ and $m_k$ to the group $G_1$ of units. The additive
structure of $Z$(.) mod $m_k$ was analysed, and extended to positive
integers by considering the principle values (naturals) of residues,
starting with $k$=3 ($Z_{30}$). Units group ~$G_1(k)$, and the additive
properties of the Boolean lattice $BL$ of idempotents of $Z_{m_k}$(.)
play an essential role.

The lower semilattice of $BL$ is multiplicative, since the $meet$ glb($a,b$)
of two idempotents is their product. The additive properties of $BL$ were
analysed, regarding the $join$ lub($a,$b) in the upper semilattice. Although
$BL$ is not closed under (+) mod $m_k$, this yields the next main results :

Lem~\ref{lem:compl}:
  ~Each complementary pair of idempotents in $Z_{m_k}$(.) sums to 1 mod $m_k$

Cor~\ref{cor:G1=G2}:
   ~Congruent max cycles $2^* \cong (2+h)^*$ ~in~ $G_2 \cong G_1$ , with $h^2 \e h=m_k/2$

Lem~\ref{lem:e+1}:
   ~Each even [odd] idempotent $e^2 \e e$ ~has~ $e$+1 in $G_1$ ~[in $G_2$]

Thm~\ref{thm:GR}: ~Each residue $2n$ mod $m_k$ is a sum of two units :
           {\it Goldbach for Residues} ~ $GR(k)$

Consider principle value set $G(k)$ of units mod $m_k$ in group $G_1(k)$.

Eqn(4):
    Restrict $2n \in G(k)+G(k)$ to prime sums $2p_{k+1} \leq 2n < (p_{k+1})^2$.

Lem~\ref{lem:GRepi}:~~Epimorphism ~$G_1(k+1) \rightarrow G_1(k)$ mod $m_k$.

Thm~\ref{thm:GC}: Goldbach Conjecture $GC$ holds, via a proof by finite reduction:\\ \hspace*{.5cm}
  $incomplete~S_0(k+1) \longrightarrow incomplete~S_0(k)$ and contradiction to complete $S_0(3)$.

\newpage
\section*{References}

\begin{enumerate}
\item T.Apostol: "{\it Introduction to Analytical Number Theory}"
  thm 10.4-6, Springer Verlag 1976.
\item N.Benschop: "The semigroup of multiplication mod $p^k$,
  an extension of Fermat's Small Theorem, and its additive structure",
 {\it Semigroups and Applications} p7, Prague, July'96.
\item N.Benschop: "Additive structure of the Group of units mod $p^k$,
 with Core and Carry concepts for extension to integers", Acta Mathematica
 Univ. Bratislava (nov.2005) \\ http://pc2.iam.fmph.uniba.sk/amuc/\_vol74n2.html
  (pp169-184, incl. direct FLT proof)
\item G.Birhoff, T.Bartee: "{\it Modern Applied Algebra}", McGraw-Hill, 1970.
\item A.Clifford, G.Preston: "{\it The Algebraic Theory of Semigroups}",~
   Vol.I, \\AMS survey \#7, p130-135, 1961.
\item S.Schwarz: "The Role of Semigroups in the Elementary Theory of
   Numbers", \\{\it Math.Slovaca} V31, N4, pp369-395, 1981.
\item K.Chandrasekharan: "{\it Introduction to Analytic Number Theory}"
   (Ch.7 - Thm 4), \\Springer Verlag, 1968.
\item N.Benschop:  "Powersums representing residues mod $p^k$, from Fermat to Waring",
   \\ {\it Computers and Mathematics, with Applications}, V39 (2000) N7-8 pp253-261.
\end{enumerate}

\end{document}